\documentclass[12pt]{article}
\usepackage[colorlinks]{hyperref}
\usepackage{color}
\usepackage{graphicx}
\usepackage{graphics}
\usepackage{makeidx}
\usepackage{showidx}
\usepackage{latexsym}
\usepackage{amssymb}
\usepackage{verbatim}
\usepackage{amsmath}
\usepackage{amsthm}
\usepackage{amsfonts}
\usepackage{amssymb,amsmath}
\usepackage{latexsym,amsthm,amscd}
\usepackage[all]{xy}

\newtheorem{prop}{Proposition}[section]
\newtheorem{defn}{Definition}[section]
\newcounter{alphthm}
\newtheorem{propriete}[alphthm]{Theorem}

\newtheorem{ex}{Example}[section]
\newtheorem{thm}{Theorem}
\newtheorem{lem}{Lemma}[section]
\newcommand{\be}{\begin{equation}}
\newcommand{\ee}{\end{equation}}
\newcommand{\ben}{\begin{enumerate}}
\newcommand{\een}{\end{enumerate}}
\newcommand{\beq}{\begin{eqnarray}}
\newcommand{\eeq}{\end{eqnarray}}
\newcommand{\beqn}{\begin{eqnarray*}}
\newcommand{\eeqn}{\end{eqnarray*}}

\newcommand{\bpf}{\begin{proof}}
\newcommand{\epf}{\end{proof}}
\newcommand{\bl}{\begin{lem}}
\newcommand{\el}{\end{lem}}
\newcommand{\bp}{\begin{prop}}
\newcommand{\ep}{\end{prop}}
\newcommand{\bd}{\begin{defn}}
\newcommand{\ed}{\end{defn}}
\newcommand{\bt}{\begin{thm}}
\newcommand{\et}{\end{thm}}
\newcommand{\bg}{\begin{ex}}
\newcommand{\eg}{\end{ex}}
\newcommand{\R}{I\!\! R}
\def\nn{\nonumber}
\newcommand\bpr{\begin{prop}}
\newcommand\epr{\end{prop}}

\title{Time optimal control in presence of moving obstacles for a Dubins airplane}
\author{Z. Fathi$^1$,  B. Bidabad$^2$\thanks{The corresponding author, bidabad@aut.ac.ir  and behroz.bidabad@math.univ-toulouse.fr.}, M. Najafpour$^3$ }
\date{\small $^{1,2}$ 
Amirkabir University of Technology (Tehran Polytechnic) Iran.\\
 $^2$Institut de Math\'{e}matique de Toulouse, Universit\'{e} Paul Sabatier, France.\\
 $^3$
University of Qu\'{e}bec at Montreal, Canada.}
\begin{document}
\maketitle
\noindent
\begin{abstract}
\noindent
In the present work, control of time-optimal trajectory for a Dubins airplane in presence of moving and fixed obstacles is obtained.
  We show that for a Dubins airplane with an initial position, the control variable can be obtained using the exact penalty function method so that the airplane reaches the end position in the shortest time in the presence of obstacles.

\end{abstract}
\vspace{.5cm}
\section{Introduction}
Optimal time control problems are widely used in engineering, particularly in robotics and aerospace.The purpose of these studies is to determine a control function that minimizes the objective function to minimize the time, trajectory or cost of a trip.

  The exact penalty methods for solving constrained optimization problems rely on the construction of a function whose unconstrained minimization points are also a solution of the constrained problem.

  The problem of Zermelo navigation or briefly ZN in Riemannian geometry is studied in detail in \cite{BRS}. In \cite{LCKTC}, the optimal time problem with fixed and moving obstacles is solved using an exact penalty function method. In a recent joint work, one of the present authors studied the time optimal trajectories of an object pursuing a moving target, without limited control in a non-obstacle space, in the context of the ZN problem, see \cite{BR}.

  A Dubins car is a simple model of a differential robot, with a constant unit speed and a minimal radius of rotation to the left or right or, equivalently, a maximum curvature equal to one.
The robot moves forward only, with the prefixed initial and terminal orientations which are prescribed with the tangent vectors to the path.

In 1957, L. E. Dubins using geometrical arguments showed that any optimal path consists of the curves with maximum curvature and straight line segments. This result was later shown using the Pontryagin's Maximum Principle Method, cf. \cite{D}.

In 2007,  H. Chitsaz and S. M. LaValle, have considered a Dubins airplane as a Dubins car having altitude. The turning angle of Dubins airplane is considered to be the turning angle of its image in the plane. Consequently, the minimum rotating radius to the left and right for Dubins airplane in space is considered equal to one.  It only flies forward and the system has independent bounded control over the altitude velocity as well as the turning rate in the plane.  Using the Pontryagin Maximum Principle, they   characterized the time-optimal trajectories for the Dubins airplane. These paths are composed of turns with minimum radius, straight line segments, and pieces of planar circular arcs, cf. \cite{CL}.

In 2010, one of the present authors in a joint work has studied a geometric approach of Dubins airplane, using Pontryagin's Maximum Principle Method and showed that its time optimal paths are geodesics of certain Finsler metric, cf. \cite{BS}.

In the present work, the time-optimal path for a Dubins airplane in presence of $n$ moving obstacles with known trajectories, from some starting point to some final point  is investigated.  To simplify the calculations, we suppose here that the number of obstacles is $n=2$. The general case is formulated in the same way.
The airplane should arrive at the final point in the shortest time without being in conflict with obstacles.
 To this aim, a control parametrization technique together with the time scaling transform is used and the problem is transformed into a sequence of optimal parameters selection problems with continuous inequality constraints and initial and terminal states constraints.
  In this problem, an airplane starts from an initial position $ (x_0, y_0, z_0) $ to reach a final position $ (x_f, y_f, z_f) $, where the trajectory of the two obstacles is previously known to the aircraft.

 Here the control variable is a $2$-tuple vector where the components are steering angle and run up movement in direction of $z$-axis. In fact, the steering angle and the moving up movement are controlled.

 The objective is to find a control such that the airplane reaches the final position in the shortest time. To solve this problem, an exact penalty function method introduced in \cite{LLTW}, \cite{YTZB} and
 \cite{YTZB2}, is used to construct a constraint violation function.
  Then, a control such as the airplane reaches the final position in shortest time is found.
  Finally, using the Matlab software, several real examples are illustrated and the effectiveness of the proposed method is illustrated.

   \section{Preliminaries and conventions}
   \subsection{Optimal control}
   In most cases, the behavior of a control system is identified by a set of differential equations that define the relationship between the input and output data.
      As a general rule, the behavior of the control system is described by an ordinary differential equation in the form of state space,
    \begin{align} \label{csi}
    w'(t)=f(w(t),{u}(t),t),
    \end{align}
     where, the $m$-tuple vector ${u}(t)$ represents the control variables at the time instant $t$.
     The phase space of a dynamic system is a space in which each possible state of the system corresponds to a single point in the phase space.
 It can normally be assumed that the phase space is an $n$-dimensional smooth manifold.
\subsection{Dubins airplane and problem statement}
In the present work, the time-optimal trajectory for a Dubins airplane in the presence of movable obstacles of known trajectories is studied.
Unlike other works, in this problem, we use an exact penalty function method and get the time-optimal trajectory for the Dubins aircraft. Moreover, we consider the presence of the $n=2$ number of moving obstacles with known trajectories. The problem for $n>2$ fix or moving obstacles is a simple extension of these computations.
In order to solve this problem, we assume that there are three moving objects, a Dubins plane and two  moving obstacles.
  The Dubins airplane is autonomously controllable and assumed to be faster than moving obstacles.
Let $w(t)=[x(t),y(t),z(t)]^T$ represents the position of the Dubins aircraft, where the parameter $ z (t) $ represents the altitude at the time $t$ in the system.
Therefore, $w'(0)$, its time derivative at the starting point could be considered as the direction tangent to the path at the starting point. For example, $w'(0)$ could be considered as the direction of the takeoff piste at the airport.
Hence its control system
 could be modeled by
\begin{align*}
 w'=(V_{xy} \cos \theta(t),V_{xy}\sin \theta(t), h'(t)),
 \end{align*}
where $V_{xy}$, is the speed in $xy$-plane, $h'(t)$ the speed in the direction of $z$-axis and $ \theta(t)$ is the angle between $x$-axis and the airplane line of sight axis in $ xy$-plane at the time t. If $ V $ denotes the speed of the airplane then we have
\begin{align}\label{zz}
 V_{xy}^2(t)+h'(t)^2=\mathbf{V}(t),  \:\:\:\: \vert u(t)\vert \leq U,  \:\:\: \forall t \geq 0,
 \end{align}
 where $u(t)$ is the 2-tuples control vector of the airplane, determined by the two variables $\theta(t)$ and $h'(t)$.
  The control vector is subject to a magnitude constraint given by (\ref{zz}).
 Let us denote by $w_1(t)=[x_1(t),y_1(t),z_1(t)]^T $  and $ w_2(t)=[x_2(t),y_2(t),z_2(t)]^T $, $\forall t \geq0$ the trajectories of the two obstacles.
 The Dubins airplane tends to fly from the starting position point $w_0$ with the initial direction $w'_0$ to the final point $w_1$, where there are two moving obstacles in its trajectory.
  The airplane should arrive at the point $w_1$ without any conflict with the moving obstacles. Assume that it arrives at the final position at the instant $T$.
If $R_1$ and $ R_2$ are the safety radiuses of the two moving obstacles and $ R $ is the safety radius for the Dubins airplane. Then the distance of the airplane from the moving obstacles should satisfies
\begin{align*}
\sqrt{[x(t)-x_i(t)]^2+[y(t)-y_i(t)]^2+[z(t)-z_i(t)]^2}\geq max\lbrace R,R_i\rbrace, \quad i=1,2.
\end{align*}
We consider here the problem of the time-optimal control, which is mathematically formulated as follows
\begin{align}\label{dap}
\displaystyle\left \{ \begin{array}{llll}
\mathrm{min}  T, \:\:\: \: \vert u(t)\vert \leq U, \quad 0\leq t\leq T,\\
w'(t)=f(w(t),u(t),t), \:\:\: w(0)=w_0, \: \:  w'(0)=w'_0,\\
\sqrt{[x(t)-x_i(t)]^2+[y(t)-y_i(t)]^2+[z(t)-z_i(t)]^2}\geq max\lbrace R,R_i\rbrace , \\
w(T)=w_1, \qquad i=1,2.
\end{array} \right.
\end{align}
We state the problem as follows. Consider a system of dynamic equation $w'(t)=f(w(t),u(t),t)$, where $ 0\leq t\leq T $ and $ w(0)=w_0 $. One of the classic problems of optimal control theory is to find a $ u(t) $ control function that minimizes the following function,
 $$\displaystyle J(u)=\psi (w(T))+\int_0^T \mathcal{L}(w(t),u(t),t)dt,$$
 where $\psi$ is the cost function and is continuously differentiable  with respect to $ w $,  $\mathcal{L}(w(t),u(t),t)$ is the Lagrangian and is continuously differentiable with respect to all arguments.
  Recall that a Lagrangian on a manifold (or a phase space) $M$ is a mapping $ \mathcal{L}:TM\longrightarrow \mathbb{R} $ which is smooth on $ TM_0 $.
 Theoretically, the principle of Pontryagin Maximum (or briefly PMP) states a solution to the problem above. Before presenting the maximum principle of Pontryagin on this work, we must remember the notion of Hamiltonian.
In general, Hamiltonian and Lagrangian mechanics are two formalisms of classical mechanics. In Hamiltonian mechanics, the trajectory of a moving particle is found without paying attention to the forces and geometry of the dynamic system.
In this way, a symplectic manifold $(M,g)$ is considered as a phase space of the dynamical system, and any real smooth function of this manifold is called a Hamiltonian.
 In physics, it is called the energy of a system.
  In Newtonian mechanics, it suffices to consider the following function as a Hamiltonian function for a fixed Lagrangian, see \cite{RMi}.
\bd
 Let $w(t)$ be the trajectory of a particle in a system with the momentum $p(t)$ and the control function ${u}(t)$. The Hamiltonian of this system with respect to the Lagrangian ${\cal L}$ is defined by:
$${\cal  H}(w(t),p(t),{u}(t))=< f(w(t),{u}(t),t),p(t)>+{\cal L}(w(t),u(t),t),$$
where $ < , >  $ is the inner product.
\ed
 The following theorem is well known
\begin{propriete} \cite[p.49]{Ev}
Let $u^*(t)$ be an optimal solution for a Hamiltonian system and $w^*(t)$  the corresponding trajectory. There exists a function
 $p^*:[0,T]\to\R^n$ such that:
\\1) we have the following Hamilton equations:
\begin{align*}
\displaystyle
\dfrac{dw^*}{dt}&=\dfrac{\partial \mathcal{H}}{\partial p}\\
\dfrac{dp^*}{dt}&=-\dfrac{\partial \mathcal{H}}{\partial w}.
\end{align*}
 2) if $\mathcal{A} = \lbrace u(t):[0,\infty) \longrightarrow \mathbb{R}^n \vert u \: is \: measurable\}$ then the conservation of energy is given by
\begin{align*}
\displaystyle \mathcal{H}(w^*(t),p^*(t),u^*(t))=\max_{u(t)\in \mathcal{A}}\mathcal{H}(w^*(t),p^*(t),u(t)).
\end{align*}
 3) the conservation of energy: $\mathcal{H}(w^*(t),p^*(t),u^*(t))$ is constant.
\\
 4) the terminal condition:
\begin{align*}
u^*(T)=\nabla \psi (w^*(T)),
\end{align*}
where, $\nabla \psi (w^*(T))$ is the gradient of cost function at the end point.
\end{propriete}
\section{Control parametrization}
It happens in the theory of optimal control that we have some path constraints and the control variables have constraints as well and that we can not obtain useful information by taking differential of the Hamiltonian function. In other words these constraints do not help much in finding $u$ or its control structures, unless in specific problems.

In present work the problem is a nonlinear optimal control problem subject to the continuous inequality constraints.  This problem is hard for solving by the classical optimal control theory.
In addition, there is an inequality constraint in each point at any given time, which implies that there is an infinite number of constraints.
\subsection{Approximate control to switch at each characteristic time}To solve this problem we shall apply the control parametrization time scaling transformation \cite{TKLJ}.
Hence the control parametrization \cite{LLTKW} is achieved as follows.
Let $t_1 $,...,$ t_p$ be the switching times, where the airplane changes its trajectory at $ t_i $ for $ i=1,...,p $. We shall employ the time scaling transform, to map these switching times into fixed time points $\frac{k}{p}$,  $ k = 1,2, . . . , p -1$, on a new time horizon $[0,1]$, see  \cite{LRTR} for more details.
Using a piecewise constant function the control function is approximated as follows
\begin{align*}
u_p(t)=\sum_{i=1}^p \chi_{_{[t_{i-1},t_i]}}(t) \sigma_i,
\end{align*}
where,
$ \sigma_{_i} $ is a 2-tuple vector and  $ t_i\leq t_{i+1} $ for $ i=1,..,p-1 $ and $ \chi_{_I} $ is the characteristic function defined for the interval $ I $  by
\begin{align*}
\chi_{_I}(t)=\displaystyle\left \{ \begin{array}{ll}1,\:\: \: \:\: \:\:\:\:\: if\:\: t\in I\\ \displaystyle 0,\:\:\:\:\: \:\:\:\:\: if \:\: t \notin I.
\end{array} \right.
\end{align*}
Next we use a rescaling time method. Let $ \Theta=\lbrace \varrho=[\rho_{_1},\rho_{_2},...,\rho_{_p}]\in \mathbb{R}^p: \rho_{_i} \geq 0, i=1,2,...,p \rbrace $, then for any $\varrho\in \Theta$ such that
\begin{align*}
\sum_{i=1}^p \frac{\rho_{_i}}{p}=T,
\end{align*}
we have a monotonic transformation from the time $ t \in [0,T] $ to a new time scale
$ s \in [0,1] $ by
\begin{align}\label{1}
 v^p(s):=\frac{dt(s)}{ds}=\sum_{k=1}^p \rho_{_k} \chi_{[\frac{k-1}{p},\frac{k}{p}]}(s), \quad s\in [0,1],
\end{align}
where,
$ t(0)=0 $.
Integrating (\ref{1}) and using the initial condition leads to
\begin{align}
t(s)=\sum_{i=1}^{k-1}\frac{\rho_i}{p}+\frac{\rho_{_k} (ps-k+1)}{p},\quad (k=1,..,p),
\end{align}
where $ s\in [\frac{k-1}{p},\frac{k}{p}] $, obviously $t(1)=\sum_{i=1}^p \frac{\rho_{_i}}{p}=T$.
Therefore, the relative subintervals are with the same size. After this rescaling, by using (\ref{1}) the control system (\ref{csi}) becomes
\begin{align}\label{I}
\frac{dw}{ds}=\frac{dw}{dt}\times \frac{dt}{ds}=v^p(s)f(w(t(s)),u(t(s)),t(s)).
\end{align}
We denote the re-scaled airplane motion by
\begin{align}\label{II}
H(s)=w(t(s)),\:\:\ H(s)=(H_1(s),H_2(s),H_3(s)).
\end{align}
Clearly for the time $ s $ in the intervals $ [\frac{k-1}{p},\frac{k}{p}]$, where $(k=1,...,p)$, the equation (\ref{1}) reduces to  $ v^p(s)=\rho_{_k}  $, so by means of (\ref{I}) and (\ref{II}) we have the following system
\begin{align*}
\displaystyle\left \{ \begin{array}{lll}\dfrac{dH}{ds}=\rho_{_k}  f(H(s),\sigma_k ,t(s)),  \quad s\in [\frac{k-1}{p},\frac{k}{p}],
\\ \displaystyle H(0)=(x(0),y(0),z(0)), \quad   H'(0)=\rho_1f(0),
\\ \displaystyle t(0)=0.
\end{array} \right.
\end{align*}
Therefore, our goal is to find $ \lbrace \sigma_k,\rho_{_k}  \rbrace $, such that $\displaystyle\sum_{k=1}^{p}\frac{\rho_{_k} }{p}$
is minimized subject to the path constraints.
\subsection{An exact penalty function method}
As mentioned earlier, this problem is an optimization problem that is subject to both the equality constraints
 \begin{align*}
 t(1)=T,\: H(1)=w_1,
  \end{align*}
  and the continuous inequality constraint
  \begin{align*}
  \Vert H(s)-H_i(s)\Vert \geq \max \lbrace R,R_i\rbrace.
  \end{align*}
 We use an exact penalty method introduced in \cite{LYTD} and \cite{LLTW}, to add all the constraints to the objective function, which poses a new problem of selecting the optimal parameters without constraint.
   Hence the new penalty functional is defined by
\begin{align}\label{cost func}
J_\delta(\sigma ,\rho ,\varepsilon)=\displaystyle\left \{ \begin{array}{lll}\displaystyle\sum_{k=1}^p\dfrac{\rho_{_k}}{p}, &\varepsilon=0,\: L(H(s))=0,
\\
 \displaystyle \sum_{k=1}^p\dfrac{\rho_{_k}}{p}+\varepsilon^{-\alpha}(L(H(s)))+\delta \varepsilon^\beta & \varepsilon>0,
\\
\displaystyle \infty & \varepsilon=0,\: L(H(s))\neq 0,
\end{array} \right.
\end{align}
where $\varepsilon \in [0,\overline{\varepsilon}] $, for $ \overline{\varepsilon}>0 $  is a new decision variable, which is a given upper bound and
\begin{align}
L(H(s))=&\displaystyle\sum_{i=1}^2\sum_{k=1}^p\int_{\frac{k-1}{p}}^{\frac{k}{p}}\rho_{_k} \max \left\lbrace \max \lbrace R^2,R_i^2\rbrace -\Vert H(s)-H_i(s)\Vert^2,0\right\rbrace ^2 ds   \nn \\ &
+\left\lbrace \sum_{i=1}^p\frac{\rho_{_i}}{p}-T\right\rbrace^2.
\end{align}
 Here, $ \delta>0 $  is the penalty parameter and $  \alpha$ and $  \beta$ are positive constants satisfying $ 0\leq \beta\leq \alpha$.
The idea of the present penalty function method can be interpreted as follows. During the process of minimizing the cost function (\ref{cost func}), whenever the penalty parameter $ \delta $ increases, $ \varepsilon^\beta $ should decrease. That is to say, $\varepsilon $ reduces as $ \beta $ is fixed. Thus, $ \varepsilon^{-\alpha} $ will increase and hence the constraint violation will reduce. This is equivalent to say the value of $ L(H(s)) $ must cut down. In this way, the satisfaction of the path constraint will eventually be achieved.
The path constraints are
\begin{align*}
\displaystyle\left \{ \begin{array}{lll}\Vert H(s)-H_i(s)\Vert^2 \geq max \lbrace R^2,R_i^2\rbrace,\: i=1,2;\: s\in [\frac{k-1}{p},\frac{k}{p}],\: k=1,2,...,p,
\\
 H(1)=w_1,
\\
t(1)=T.
\end{array} \right.
\end{align*}
Next, the goal is to find $(\sigma,\rho,\varepsilon)$ such that the cost function $J_{\delta}(\sigma,\rho,\varepsilon)$ is minimized subject to $\varepsilon\geq 0 $. Hence the control variable can be obtained such that the airplane reaches the final position in shortest time.
Now, the objective function of this problem is in canonical form. To solve this problem, we need the gradient formula of the objective function. As derived in the proof of Theorem $5.2.1$ in \cite{KTGK}, it is well known the gradient formulas are given in the following form.
For each $ \delta>0 $, the gradients of the cost function $ J_\delta(\sigma,\rho,\varepsilon) $ with respect to $ \sigma $ and $ \rho $ are:
\begin{align*}
\dfrac{\partial J_\delta(\sigma,\rho,\varepsilon)}{\partial \sigma}=\dfrac{\partial \psi_{_0} \big{(}\rho,\varepsilon,H(s)\big{)}}{\partial \sigma}+\int_0^1 \dfrac{\partial  \mathcal{H}_0\big{(}s,H(s),\sigma,\rho,\lambda(s)\big{)}}{\partial \sigma}ds,
\end{align*}
\begin{align*}
\dfrac{\partial J_\delta(\sigma,\rho,\varepsilon)}{\partial \rho}=\dfrac{\partial \psi_{_0} \big{(}\rho,\varepsilon,H(s)\big{)}}{\partial \rho}+\int_0^1 \dfrac{\partial \mathcal{H}_0\big{(}s,H(s),\sigma,\rho,\lambda(s)\big{)}}{\partial \rho}ds,
\end{align*}
where, $\mathcal{H}_0 \big{(}s,H(s),\sigma,\rho,\varepsilon\big{)}$ is the Hamiltonian of the cost function given by
\begin{align*}
\mathcal{H}_0 \big{(}s,H(s),\sigma,\rho,\varepsilon\big{)}= \displaystyle\sum_{i=1}^2\sum_{k=1}^p\mathcal{L}_{0,ik}\big{(}H(s),\rho\big{)}
+\lambda_{0}(s)v^p(s)f\big{(}H(s),\sigma,t(s)\big{)},
\end{align*}
  therein $\psi_{_0} \big{(}\rho,\varepsilon,H(s)\big{)}$ and $\mathcal{L}_{0,ik}\big{(}H(s),\rho\big{)}$ are defined by
\begin{align*}
&\psi_{_0} \big{(}\rho,\varepsilon, H(s)\big{)}=\displaystyle\sum_{i=1}^p\frac{\rho_k}{p}+\delta \varepsilon^\beta+\varepsilon^{-\alpha}\left\lbrace \sum_{i=1}^p\frac{\rho_{_i}}{p}-T\right\rbrace^2,\\
&\mathcal{L}_{0,ik}\big{(}H(s),\rho \big{)}=\rho_k max \left\lbrace \max \lbrace R^2,R_i^2\rbrace -\Vert H(s)-H_i(s)\Vert^2,0\right\rbrace ^2.
\end{align*}
Here, $\lambda_0(s)$ is the solution of the following co-state differential equation
\begin{align*}
(\lambda_0(s))^T=\dfrac{\partial \mathcal{H} _0\big{(}s,H(s),\sigma,\theta,\lambda(s)\big{)}}{\partial H(s)},
\end{align*}
with the boundary condition
$(\lambda_0(1))^T=\dfrac{\partial \psi_0\big{(}\rho,\varepsilon,H(s)\big{)}}{\partial H(s)}$.
\subsection{Computational results}
\begin{figure}[hth]
\begin{center}
\includegraphics[scale=.4]{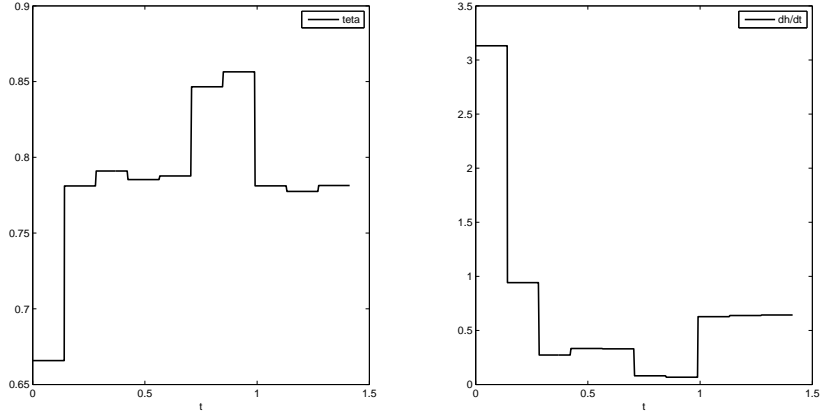}
\caption{\small Optimal control $\theta$ and $h'=\frac{dh}{dt}$ for Dubins airplane from $ (0,0,0) $ to (1,1,1), where the obstacles are far from the airplane and $V_{xy}=1\frac{m}{s}$, $T=1.4159$.}
\label{ex1+.eps}
\end{center}
\end{figure}

\begin{figure}[hth]
\begin{center}
\includegraphics[scale=.4]{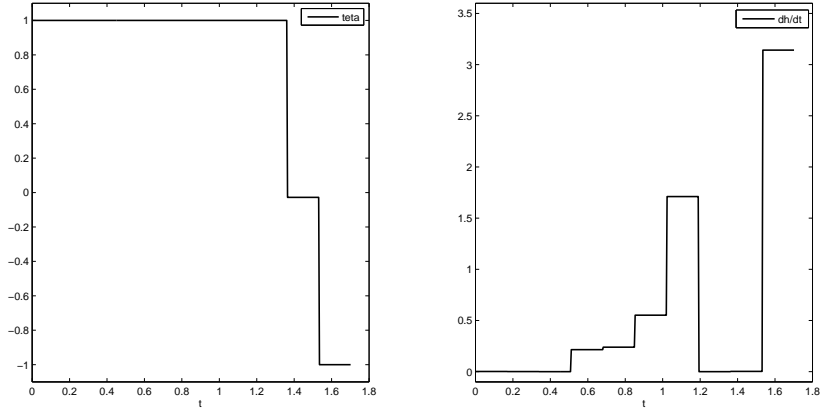}
\caption{Optimal control $ \theta $ and $ h' $ for Dubins airplane where the obstacles are on the trajectory, $ V_{xy}=1 \frac{m}{s}$, $ T=1.7058 $.\small}
\label{ex2.eps}
\end{center}
\end{figure}

\begin{figure}[hth]
\begin{center}
\includegraphics[scale=.4]{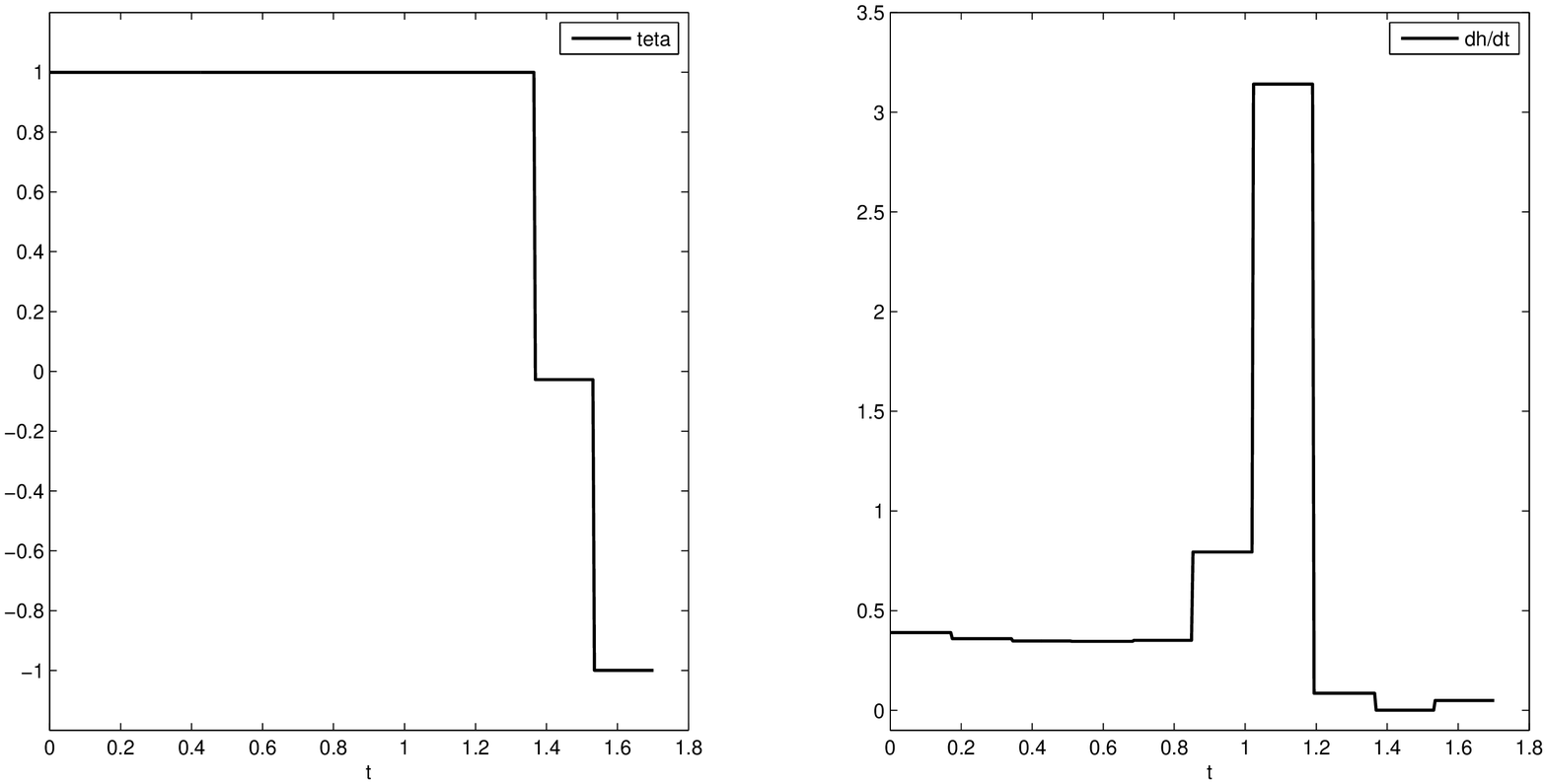}
\caption{Optimal control $\theta$ and $h'$ for Dubins airplane, $V_{xy}=1 \frac{m}{s}$, $T=1.7059$.\small}
\label{ex3.eps}
\end{center}
\end{figure}

 In this section, using Matlab software program we present three different examples. In the first example, it is assumed that the moving obstacles have not important effects on the trajectories,  since they are far from  the path of airplane. In the next examples moving obstacles are laying along the trajectory of airplane.
\bg[Non-important obstacles]\label{no}
The Dubins airplane tends to fly from the starting point $(0,0,0)$ to the final $(1,1,1)$. The two moving obstacles  $A$ and  $B$ are not near the airplane trajectory. Assume the trajectories of the two obstacles $A$ and $B$ are described as:
\begin{align*}
&A(t)=(t,t, \sin\frac{\pi t}{5}), \:\:\:\:\:\: t\geq \frac{1}{10} \\
&B(t)=\displaystyle\left \{ \begin{array}{ll}\displaystyle (\frac{1}{2}-\sqrt{\frac{3}{4}-2(t-\frac{1}{2})},t,t), \ \:\:\:\:\:  \frac{1}{10}\leq t\leq \frac{\sqrt{3}}{8}+\frac{1}{2},
\\
\displaystyle (\frac{1}{2}+\sqrt{\frac{3}{4}-2(t-\frac{1}{2})},t,t), \:\:\:\:\ \frac{\sqrt{3}}{8}+\frac{1}{2}\leq t\leq 1.
\end{array} \right.
\end{align*}
Let the motion equation of the airplane be
\begin{align*}
 w'=(V_{xy} \cos \theta(t),V_{xy}\sin \theta(t), h'(t)),
\end{align*}
where, $V_{xy}=1 \:(m/s)$ is the velocity in $xy$-plane and the radius of safety region for airplane is $ 0.2 $ and the safety region for moving obstacles is $0.1$. In this example the airplane flies from the initial point $(0,0,0)$ to the final point $(1,1,1)$. We obtain the control variables $h'$ and $\theta(t)$ such that the airplane reaches $(1,1,1)$ in shortest time without encountering the obstacles $A$ and $B$. The optimal control $h'$ and $\theta$ are shown in Fig \ref{ex1+.eps}, drawn using Matlab software program.
In this example $ \alpha=1 $, $\delta=50  $ and $ \beta=1 $.
\eg
In the next example we consider a moving obstacle with a prefix path.
\bg Let the trajectory of one of the obstacles be a straight line between starting and final points.
With the hypothesis of the example (\ref{no}) on Daubins airplane and the moving obstacle $A$,  we have put $B(t)=(t,t,t),$  for $t\geq \: \frac{1}{10}$, $V_{xy}=1 \frac{m}{s}$ and the safety region equal $0.1$. Then the
optimal control $h'$ and $\theta$ are shown in Fig \ref{ex2.eps}, drawn using Matlab software program.
In this example $\alpha=1$, $\delta=50$ and $\beta=1$.
\eg
\bg[General example]
Let Dubins airplane fly from $ (0,0,0) $ to $ (1,1,1) $ and the trajectories of the two obstacles $ A $ and $ B $ be described
respectively
\begin{align*}
A(t)=(t,t,t), \quad t\geq \frac{1}{10}\:\:\: B(t)=(t,t,2-t),
\end{align*}
and $ A $ flies from $ (\frac{1}{10},\frac{1}{10},\frac{1}{10}) $ to $ (1,1,1) $ and $ B $ flies from $ (0,0,2) $ to $ (1,1,1) $. Let the safety radius be $0.1$. The
optimal control $ h' $ and $ \theta $ are shown in Fig \ref{ex3.eps}, drawn using Matlab software program.
In this example $\alpha=1$, $\delta=10$ and $\beta=1$.
\eg
\newpage
\section{Conclusions}
The problem of the time-optimal paths for a Dubins plane is studied in \cite{CL} and the geometry of its movement in a non-obstacle space is studied in \cite{TMRB}.

 In this paper, we have discussed an effective calculation method, namely an Exact Penalty Function Method, to determine the optimal control over time for a Dubins aircraft, from one starting point to a final point, in the presence of fixed or moving obstacles, with known trajectories.

 A numerical method is applied to three scenarios for the airplane in $3$D spaces and the time-optimal controls are established.
The results obtained clearly demonstrate the effectiveness of the proposed method.


  This work is an extension of \cite{BS, CL} and \cite{LCKTC} in certain senses.
   It can be shown that the metric associated with time-optimal trajectories  in the presence of moving obstacles is a Finsler metric. This will appear in future work.



Faculty of Mathematics and Computer Sciences,
Amirkabir University of Technology (Tehran Polytechnic),
424 Hafez Ave. 15914 Tehran, Iran.\\
 bidabad@aut.ac.ir;  z.fathi@aut.ac.ir \\
 Institut de Math\'{e}matique de Toulouse, Universit\'{e} Paul Sabatier,
F-31062 Toulouse cedex 9, France.\\
behroz.bidabad@math.univ-toulouse.fr\\
Dept. of Mathematics,
Universit\'{e} du Qu\'{e}bec   \`{a} Montr\'{e}al,\\
 405 Rue St-Catherine Est, Montr\'{e}al, QC H2L 2C4, Canada.
\\
najafpour\underline{\hspace{0.2 cm}}ghazvini.mehrdad@courrier.uqam.ca.
\end{document}